\documentclass[12pt,reqno]{amsart}
\usepackage[margin=0.5in]{geometry}
\usepackage[alphabetic]{amsrefs}
\usepackage{tabularx}
\usepackage{yfonts}
\newcommand{\R}{\mathbb{R}}

\renewcommand{\epsilon}{\varepsilon}

\usepackage[alphabetic]{amsrefs}
\usepackage{amsmath,amssymb,amsfonts,amsthm,enumerate}
\usepackage{mathrsfs}
\usepackage{mathtools} 
\usepackage{url}
\usepackage{graphicx,epstopdf,color}

\newtheorem{theorem}{Theorem}[section]

\newtheorem{corollary}[theorem]{Corollary}

\newtheorem{remark}[theorem]{Remark}
\newtheorem{conjecture}[theorem]{Conjecture}

\renewcommand{\leq}{\leqslant}
\renewcommand{\le}{\leqslant}
\renewcommand{\geq}{\geqslant}
\renewcommand{\ge}{\geqslant}

\newcommand{\abs}[1]{\left\lvert#1\right\rvert}

\begin{document}

\author{Serena Dipierro}
\author{Pietro Miraglio}
\author{Enrico Valdinoci}

\thanks{{\em Serena Dipierro}: Department of Mathematics and Statistics,
University of Western Australia,
35 Stirling Highway,
Crawley WA 6009, Australia. {\tt serena.dipierro@uwa.edu.au} }

\thanks{{\em Pietro Miraglio}:
Dipartimento di Matematica, Universit\`a degli studi di Milano,
Via Saldini 50, 20133 Milan, Italy, and
Universitat Polit\`ecnica de Catalunya,
Departament de Matem\`atiques,
Avinguda Diagonal 647, 08028 Barcelona, Spain.
{\tt pietro.miraglio@unimi.it}}

\thanks{{\em Enrico Valdinoci}: Department of Mathematics and Statistics,
University of Western Australia,
35 Stirling Highway,
Crawley WA 6009, Australia. {\tt enrico.valdinoci@uwa.edu.au} }

\title[Symmetry result in water waves models]{Symmetry results for the solutions \\ of
a partial differential equation \\ arising in water waves}

\begin{abstract}
This paper recalls some classical motivations in fluid dynamics leading
to a partial differential equation which is prescribed on a domain
whose boundary possesses two connected components, one endowed with a Dirichlet
datum, and the other endowed with a Neumann datum.

The problem can also be reformulated as a nonlocal problem on the component
endowed with the Dirichlet datum. A series of recent symmetry results
are presented and compared with the existing literature.
\end{abstract}

\maketitle

\section{Introduction}
In this paper we present some recent results related to the partial differential
equation
\begin{align}
\label{lvXsis}
\begin{cases}
\mathrm{div}(y^a\nabla v)=0 \qquad \text{for}\,\,&x\in\R^n,y\in(0,1),\\
v_y(x,1)=0 \qquad &x\in\R^n,y=1,\\
v(x,0)=u(x) \qquad &x\in\R^n,y=0,\\
-\displaystyle\lim_{y\rightarrow 0}y^av_y=f(v)\qquad &x\in\R^n,y=0,
\end{cases}
\end{align}
with~$a\in(-1,1)$. The problem in~\eqref{lvXsis}
is related to a water waves model and, in a suitable limit, it
recovers a fractional Laplace operator. More precisely,
a solution~$v$ of~\eqref{lvXsis} can be related to its trace~$u$
by a nonlocal equation of the type
\begin{equation}\label{7uhn88} \mathcal{L}_au=f(u) \qquad \text{in}\,\,\R^n,\end{equation}
for a suitable linear operator~$\mathcal{L}_a$. The operator~$\mathcal{L}_a$
can be written in Fourier modes and will present different asymptotic behaviours for small
and large frequencies, making the problem particularly interesting.

One of the main questions that we address is under which conditions
the bounded and monotone solutions of~\eqref{7uhn88}
are necessarily
one-dimensional --- that is, as a counterpart,
solutions of~\eqref{lvXsis} that are monotone in one of the $x$-variables
are necessarily functions only of~$x_1$ and~$y$, up to a rotation.\medskip

In Section~\ref{model_sec} we recall some basic fluid dynamics
motivations to give an elementary but exhaustive description of the
problem in~\eqref{lvXsis} in terms of classical physics. Then, in Section~\ref{SxSYRESU},
we focus on the mathematics relative to~\eqref{lvXsis} and~\eqref{7uhn88}, discussing
symmetry results in the light of a classical conjecture by Ennio De Giorgi.

\section{Physical considerations}\label{model_sec}

\begin{figure}
    \centering
    \includegraphics[width=12cm]{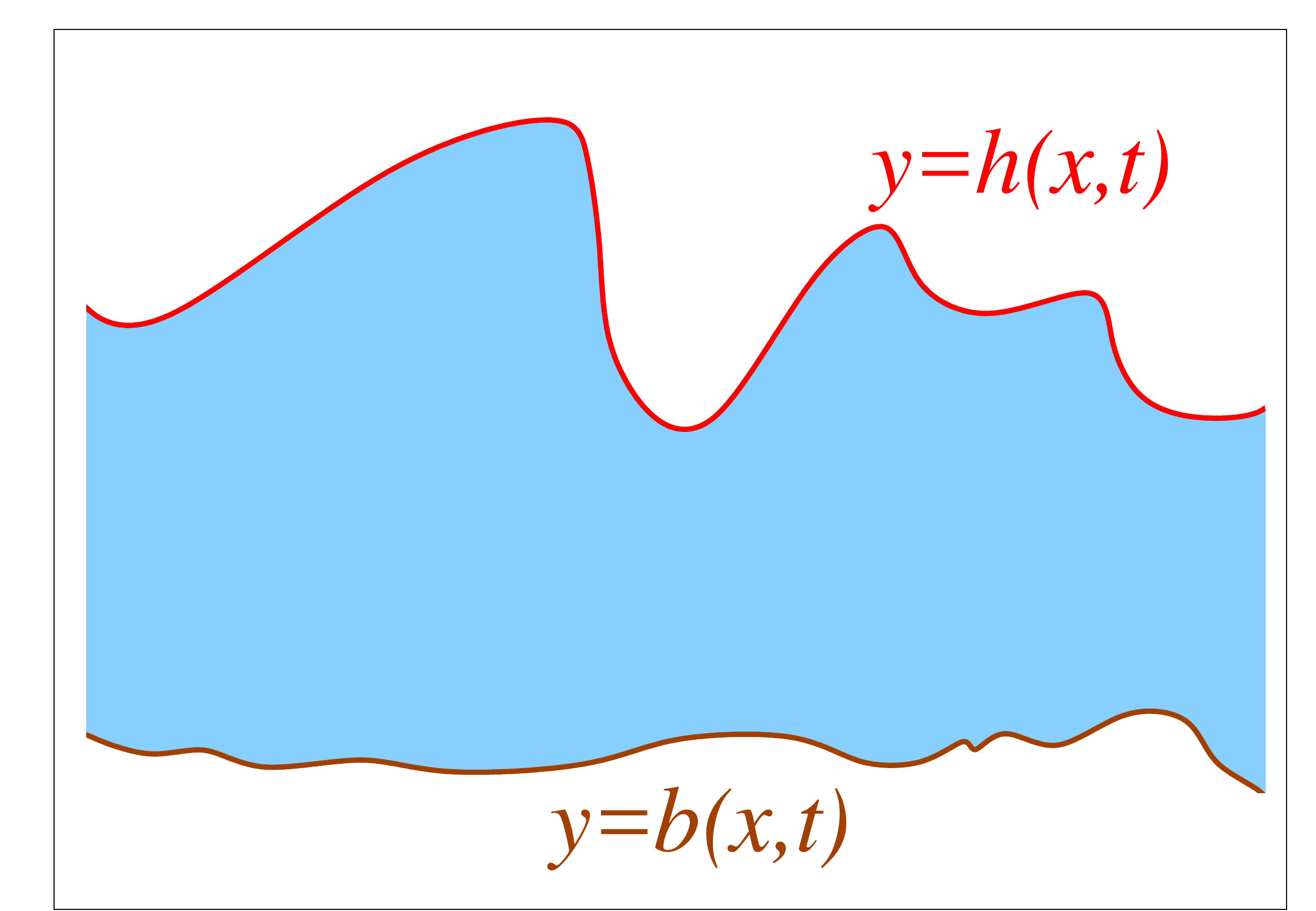}
    \caption{\em {{The domain~$\Omega(t)$ in~\eqref{MARE}.}}}
    \label{C1}
\end{figure}

In this section, we give a detailed motivation for the problem in~\eqref{lvXsis}
arising from a fluid dynamics model.
To this end, we consider a possible physical description of an irrotational and inviscid fluid
(the ``ocean'') in~$\R^{n+1}$, though we commonly take~$n=2$
in the ``real world''. The position of a fluid particle at time~$t$
will be denoted by~$X(t)=(x(t),y(t))\in\R^n\times\R$.
We suppose that, at time~$t$, the region occupied by the ocean
lies above the graph of a function~$b(\cdot,t)$ (the ``bottom
of the ocean'') and below
the graph of a function~$h(\cdot,t)$ (the ``surface of the ocean'').
Therefore, in this model, the ocean can be described by
the time-dependent domain
\begin{equation}\label{MARE}\Omega(t):=\{ (x,y)\in\R^n\times\R {\mbox{ s.t. }}
b(x,t)\le y\le h(x,t)\},\end{equation}
see Figure~\ref{C1}.

Given a point~$X\in\Omega(t)$, we denote by~$v(X,t)$
the velocity of the fluid particle at~$X$ at time~$t$.
We denote by~$\Phi^t(X)$ the evolution
produced by the vector field~$v$ at time~$t$ starting at the point~$X$
at time zero,
that is the solution
of the initial value problem
\begin{equation}\label{QUI} \begin{dcases}
&\displaystyle\frac{d}{dt} \Phi^t(X)= v(\Phi^t(X),t)\qquad{\mbox{ for (small) }}t>0,
\\ &\Phi^0(X)=X.\end{dcases}\end{equation}

We suppose that the density of the water is described by a positive
function~$\rho=\rho(X,t)$. Then, the mass of
the fluid lying in a region~$\widetilde\Omega\subset\R^{n+1}$ at time~$t$
is described by the quantity
\begin{equation}\label{MASS} \int_{\widetilde\Omega} \rho(X,t)\,dX.\end{equation}
The rate at which a fluid mass enters in~$\widetilde\Omega$
through an infinitesimal
portion of~$\partial\widetilde\Omega$ in the vicinity of a point~$X\in\partial\widetilde\Omega$ is given by the density times the velocity
at~$X$ in the direction of the inner normal of~$\partial\widetilde\Omega$
at~$X$. That is, if~$\nu(X)$ denotes the exterior normal
of~$\partial\widetilde\Omega$ at~$X$, we find that the
rate at which a fluid mass enters in~$\widetilde\Omega$ is given by
$$ -\int_{ \partial\widetilde\Omega } \rho(X,t)\,v(X,t)\cdot\nu(X)\,d{\mathcal{H}}^n(X).$$
Comparing with~\eqref{MASS}, and using the Divergence Theorem,
this leads to
$$ \int_{\widetilde\Omega} \partial_t\rho(X,t)\,dX=
\frac{d}{dt} \int_{\widetilde\Omega} \rho(X,t)\,dX=
-\int_{ \partial\widetilde\Omega } \rho(X,t)\,v(X,t)\cdot\nu(X)\,d{\mathcal{H}}^n(X)
=
-\int_{\widetilde\Omega } {\rm div}_X\big(\rho(X,t)\,v(X,t)\big)\,dX.$$
{F}rom this, since the volume region~$\widetilde\Omega$ is arbitrary,
we obtain the ``mass conservation law'' (also known as
``continuity equation'') given by
\begin{equation}\label{MCL}
\partial_t \rho(X,t)+{\rm div}_X\big(\rho(X,t)\,v(X,t)\big)=0\qquad{\mbox{ in }}\;\Omega(t).
\end{equation}

Let us now analyze the conditions occurring at the bottom and at the surface of the fluid.
At the bottom, we assume that the fluid cannot penetrate inside the ground,
hence its velocity is tangent to the seabed. Recalling the notation in~\eqref{MARE},
we have that~$v$ needs to be orthogonal to the normal direction of the graph of~$b$, and thus,
using
the notation~$X=(x,y)\in\R^n\times\R$,
\begin{equation}\label{FLU2}
v(X,t)\cdot \big(\nabla_x b(x,t), -1\big)=0\qquad{\mbox{ if }}\;y=b(x,t).
\end{equation}
We can therefore collect the results in~\eqref{MCL} and~\eqref{FLU2} by writing
\begin{equation}\label{FLU4a}
\begin{cases}
\partial_t \rho(X,t)+{\rm div}_X\big(\rho(X,t)\,v(X,t)\big)
=0&\qquad{\mbox{ in }}\;\Omega(t),\\
v(X,t)\cdot \big(\nabla_x b(x,t), -1\big)=0&\qquad{\mbox{ on }}\;\{y=b(x,t)\}.
\end{cases}
\end{equation}
{F}rom~\eqref{FLU4a} one sees that the vector field~$\rho v$
has perhaps more physical meaning than~$v$ alone, since it represents
the density speed of the flow, and it is somehow more meaningful to prescribe a bound
on~$\rho v$ rather than on~$v$ itself. For instance, the situation in which~$v$
becomes unbounded becomes physically realistic if~$\rho v$ remains bounded,
since, in this case, roughly speaking, only a very negligible amount of fluid
would travel at exceptionally high speed. Therefore, though the equations
are perfectly equivalent in case of ``nice'' vector fields~$v$ and densities~$\rho$,
we prefer to write~\eqref{FLU4a} in a form which makes appear directly the quantity~$\rho v$
rather than~$v$ alone. This is done by multiplying the identity on the bottom of the ocean by the density,
to find
\begin{equation}\label{FLU4}
\begin{cases}
\partial_t \rho(X,t)+{\rm div}_X\big(\rho(X,t)\,v(X,t)\big)
=0&\qquad{\mbox{ in }}\;\Omega(t),\\
\rho(X,t)\,v(X,t)\cdot \big(\nabla_x b(x,t), -1\big)=0&\qquad{\mbox{ on }}\;\{y=b(x,t)\}.
\end{cases}
\end{equation}

\begin{figure}
    \centering
    \includegraphics[width=12cm]{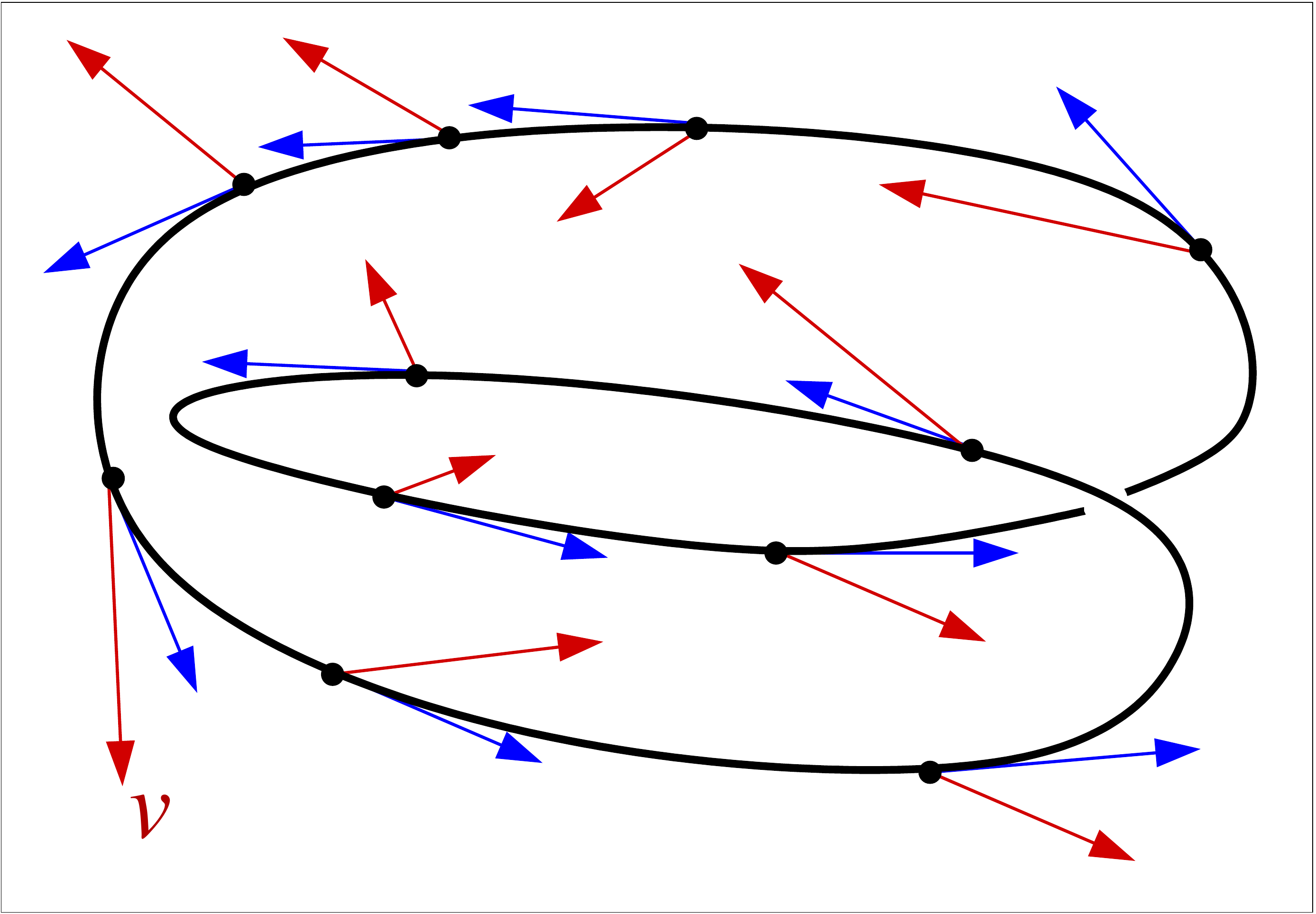}
    \caption{\em {{The velocity field~$v$ has always a positive
    component along the tangential direction of the closed curve, hence 
    it is not irrotational.}}}
    \label{C2}
\end{figure}

We also assume that the fluid particles do not ``circulate
in a cyclone way'', namely that the fluid is irrotational, see Figure~\ref{C2}.
To formalize this notion in an arbitrarily large number of dimensions in an elementary geometric
way (without using the notion of higher dimensional curls), we assume that, for every fixed time,
the integral of the velocity vector field along any closed one-dimensional curve in~$\R^n$
vanishes. As a matter of fact, it would be enough to require such a condition
along polygonal lines, and in fact it would be sufficient to require it along
triangular connections.

This irrotationality condition implies (and, in fact, 
it is equivalent to) that the velocity field admits a potential, namely that there exists
a scalar function~$u=u(X,t)$
such that
\begin{equation}\label{POTE}
v(X,t)=\nabla_Xu(X,t)
.\end{equation}
We stress that~\eqref{POTE} is a rather striking formula, since it reduces
the knowledge of a vector valued function (namely, $v$) to the knowledge of (the derivatives of)
a single scalar function. The construction of the potential~$u$ is standard, and can be
performed along the following argument: we let~$\Gamma_X$ be
the oriented segment
starting at the origin and arriving at~$X$, and we set
$$ u(X,t):=\int_{\Gamma_X} v:=\int_0^1 v(\vartheta X,t)\cdot X\,d\vartheta.$$
To prove~\eqref{POTE}, let~$j\in\{1,\dots,n\}$ and~$\delta\ne0$, to be taken arbitrarily small in what follows. We also denote by~$\Gamma_{X,\delta,j}$ the oriented segment from~$X$ to~$X+\delta e_j$.
Also, given two adjacent segments~$\Gamma_1$ and~$\Gamma_2$,
we denote by~$\Gamma_1\cup\Gamma_2$ the broken line joining the initial point of~$\Gamma_1$
to the end point of~$\Gamma_1$ (which coincides with the initial point of~$\Gamma_2$)
and that to the end point of~$\Gamma_2$. Furthermore, we denote by~$-\Gamma_1$
the segment~$\Gamma_1$ run in the opposite direction. With this notation, we have that~$\Gamma_{X+\delta e_j}\cup(-\Gamma_{X,\delta,j})\cup (-\Gamma_X)$ forms a close triangle and accordingly,
by the irrotationality condition,
\begin{eqnarray*}&&
0= \int_{ \Gamma_{X+\delta e_j}\cup(-\Gamma_{X,\delta,j})\cup (-\Gamma_X) }v =
\int_{\Gamma_{X+\delta e_j}} v-\int_{\Gamma_{X,\delta,j}}v-\int_{\Gamma_X}v\\&&\qquad=
u(X+\delta e_j,t)-\delta\int_{0}^1v(X+\vartheta \delta e_j,t)\cdot e_j\,d\vartheta - u(X,t)
.\end{eqnarray*}
Dividing by~$\delta$ and sending~$\delta\to0$, we obtain~\eqref{POTE}, as desired.

Then, inserting~\eqref{POTE} into~\eqref{FLU4}, we conclude that
\begin{equation}\label{finalmodel}
\begin{cases}
\partial_t \rho(X,t)+{\rm div}_X\big(\rho(X,t)\,\nabla_X u(X,t)\big)=0&
\qquad{\mbox{ in }}\;\Omega(t),\\
\rho(X,t)\,\nabla_xu(X,t)\cdot \nabla_x b(x,t)
-\rho(X,t)\,\partial_yu(X,t)=0&\qquad{\mbox{ on }}\;\{y=b(x,t)\}.
\end{cases}
\end{equation}
We observe that the setting in~\eqref{lvXsis} is a particular case
of that in~\eqref{finalmodel}, in which one considers the steady case
of stationary solutions (i.e.~$\rho$ does not depend on time),
with~$X=(x,y)\in\Omega=\R^n\times(0,1)$,
and~$\rho(X)=y^a$, with~$a\in(-1,1)$.

\begin{remark} {\rm
Concerning the setting in~\eqref{QUI},
we recall that in the literature
one also considers the ``streamlines''
of the fluid, described
by a parameter~$\tau\in\R$, which
are (local) solutions of the differential equation (for fixed time~$t$)
$$ \frac{d}{d\tau} X(\tau,t)=v(X(\tau,t),t).
$$
Notice that, if the velocity field~$v$ is independent of time,
we can actually identify the curve parameter~$\tau$ with the usual time~$t$ and then the streamlines describe the physical trajectories
of the fluid particle. But in general, for
velocity fields which depend on time, streamlines do not represent
the physical trajectories. Nevertheless,
streamlines
are always instantaneously tangent to the velocity field of the flow
and therefore they indicate the direction in which the fluid particle
at a given point travels in time. We maintain the distinction between streamlines
and physical trajectories of the flow, and in this note only the latter objects will be taken
into account for the main computations.
}\end{remark}

\begin{remark}{\rm We point out that in the literature
one often assumes that the fluid is ``incompressible'', that is, fixed
any reference domain~$\widetilde\Omega$,
\begin{equation*}
\frac{d}{dt}\int_{\widetilde\Omega} \rho(\Phi^t(X),t)\,dX=0.
\end{equation*}
This condition together with~\eqref{QUI} leads to
\begin{equation}\label{INC1} \partial_t\rho(\Phi^t(X),t)+\nabla_X\rho(\Phi^t(X),t)\cdot v(\Phi^t(X),t)=0,\end{equation}
or, equivalently, changing the name of the space variable
\begin{equation}\label{INC}
\partial_t\rho(X,t)+\nabla_X\rho(X,t)\cdot v(X,t)=0.\end{equation}
The incompressibility condition~\eqref{INC} may be also understood from a ``discrete analogue''
by thinking that the density~$\rho(X,t)$ of a gas formed by indistinguishable molecules at a point~$X$
at time~$t$ is measured by ``counting the number of molecules'' in the vicinity of~$X$ at time~$t$.
That is, fixing~$r>0$, the gas density could be defined as the number
of molecules lying in~$B_r(X)$ at time~$t$. If the gas is incompressible, we expect that the
number of molecules around the evolution~$\Phi^t(X)$ of~$X$ remains the same. This gives that
$$ \rho(\Phi^t(X),t)=\rho(X,0),$$
which leads to~\eqref{INC1} and~\eqref{INC}.

To appreciate the structural difference between the 
mass conservation law in~\eqref{MCL}
and the incompressibility condition in~\eqref{INC}, let us consider two examples.
In the first example, let
$$ v(X,t):=-X\qquad{\mbox{ and }}\qquad\rho(X,t):=e^{nt},$$
with $ n>0 $. In this case, the velocity field pushes all the fluid towards the origin,
preserving the mass according to~\eqref{MCL}: as a consequence, the particles
of the fluid get ``packed'' and their density increases, and the incompressibility condition~\eqref{INC}
is indeed violated.

As a second example, let us consider the case in which
$$ v(X,t):=-X\qquad{\mbox{ and }}\qquad\rho(X,t):=1.$$
In this case, the fluid elements are still pushed towards the origin, but the density remains
constant. This means that there must be a leak somewhere, from which the fluid escapes.
In this situation, the incompressibility condition in~\eqref{INC} is satisfied, but the mass is lost
and accordingly~\eqref{MCL} does not hold.

We also point out that if
the the mass conservation law in~\eqref{MCL}
and the incompressibility condition in~\eqref{INC} are both satisfied, then
\begin{eqnarray*} 0&=&
\partial_t \rho(X,t)+{\rm div}_X\big(\rho(X,t)\,v(X,t)\big)\\&=&
\partial_t \rho(X,t)+\nabla_X \rho(X,t)\cdot v(X,t)+\rho(X,t)\,{\rm div}_X\,v(X,t)\\&=&
\rho(X,t)\,{\rm div}_X\,v(X,t),\end{eqnarray*}
and, as a consequence,
\begin{equation*}
{\rm div}_X\,v(X,t)=0\qquad{\mbox{ in }}\;\Omega(t).
\end{equation*}
In this note, we will not explicitly take into account incompressibility assumptions, but merely the conservation of mass in~\eqref{MCL}.
}\end{remark}

\begin{remark}{\rm
Concerning the top surface of the fluid, in the literature it is often
assumed that fluid particles on this surface remain
there forever (i.e., there is no ``mixing effect'' between the top surface of the sea and the 
rest of the water mass). This condition, in the notation of~\eqref{MARE}
and~\eqref{QUI}, would translate into
$$ \Phi^t_2(X)=h(\Phi_1^t(X),t),$$
as long as~$X=(x,y)$ and~$y=h(x,0)$,
where~$\Phi^t(X)=\big(\Phi^t_1(X),\Phi^t_2(X)\big)\in\R^n\times\R$. Hence, in view of~\eqref{QUI},
\begin{equation*}
0=\frac{d}{dt} \Big( h(\Phi_1^t(X),t)-\Phi^t_2(X)\Big)=v(\Phi^t(X),t)\cdot
\big(\nabla_x h(\Phi^t_1(X),t),\,-1\big)+\partial_t h(\Phi_1^t(X),t).
\end{equation*}
In this note, we do not need to assume this additional no mixing condition.
}\end{remark}

\section{Symmetry results} \label{SxSYRESU}
Now, we present some results for an elliptic problem related to the
stationary case of the model introduced in Section \ref{model_sec}.
Besides assuming no dependence on time $ t $, we also consider the
simplification of a ``flat ocean'', by taking $ b(x)= H>0 $ and $ h(x)=0 $
(recall the notation in~\eqref{MARE}). This choice implies that we now consider the sea as
\begin{equation*}\Omega=\{ (x,y)\in\R^n\times\R {\mbox{ s.t. }}
0\le y\le H\},\end{equation*}
and that we are ``reversing
the vertical direction'', in order to have the ocean surface on $\{y=0\}$. This last simplification is done for pure mathematical convenience and does not affect the model.

In our setting, we can use (\ref{finalmodel}) in order to associate a velocity potential in the whole slab $ \R^2\times[0,H] $ with a given datum on the surface of the ocean. Given the values of the velocity potential on $ \{y=0\} $ and denoting
such datum by~$u$, we consider the velocity potential $ v $ in the whole slab $ \R^2\times[0,H] $ that solves	
\begin{align}\label{1st}
\begin{cases}
0={\rm div }(\varrho V)={\rm div }(\varrho\nabla v ) 
\quad&{\mbox{ in }}\R^2\times(0,H),\\
0=V_3\big|_{y=H}= v_y\big|_{y=H}\quad&{\mbox{ on }}
\R^2\times\{y=H\} ,\\
v=u\big|_{y=0}\quad&{\mbox{ on }}\R^2\times\{y=0\}.
\end{cases}
\end{align}
In relation to water waves and in view of the discussion in Section \ref{model_sec}, we
are interested in the weighted vertical velocity on the surface of the ocean. Thus,
the operator that we want to study is
\begin{equation}\label{op}
\mathcal{L}_a u(x):=-\displaystyle\lim_{y\rightarrow 0}\rho(y)v_y(x,y).
\end{equation}
When~$\varrho:=1$ and~$H\to+\infty$ (which is the case of a fluid with constant density
and an ``infinitely deep sea''), the operator $ \mathcal{L}_a $ is the square root of
the Laplacian, see e.g.~\cite{CS}. For finite values of~$H$
the operator described in~\eqref{op} is nonlocal, but also not of purely fractional type, as we are going to see.

In the following, 
we choose
\begin{equation}\label{RRHP}
\varrho(y):=y^a\end{equation}
as a density, where $a\in (-1,1)$.
We notice that, in this case, 
\begin{equation}\label{ST}
{\mbox{the limit as~$H\to+\infty$ corresponds to the $s$-th root
of the Laplacian,}}\end{equation} with~$s:=(1-a)/2$, but for a finite value of~$H$
the problem is not of purely fractional type. From now on, we normalize the domain by setting~$H:=1$. 
{F}rom a physical point of view, the choice in~\eqref{RRHP}
corresponds to the situation in which the
density of the fluid at a point depends only on the
depth, in a power-like fashion,
and it is constant in the horizontal directions.
Possibly, some of the results that we present here can be extended to the
case of a more general density $ \rho(y) $, and
we intend to investigate the possibility of this generalization in a forthcoming work.

After generalizing the physical setting $\R^2\times[0,1] $ to the
mathematically
interesting case $\R^n\times[0,1]$ --- with coordinates $x\in\R^n$ and~$ y\in[0,1]$ --- the
extension problem in~(\ref{1st}) reads
\begin{align}
\label{lvsis}
\begin{cases}
\mathrm{div}(y^a\nabla v)=0 \qquad &\text{in}\,\,\R^n\times (0,1),\\
v_y(x,1)=0 \qquad &\text{on}\,\,\R^n\times \{y=1\},\\
v(x,0)=u(x) \qquad &\text{on}\,\,\R^n\times \{y=0\}.
\end{cases}
\end{align}
Therefore, in light of~\eqref{RRHP},
the Dirichlet to Neumann operator $\mathcal{L}_a$ in~\eqref{op}
is given by 
\begin{equation*}
\mathcal{L}_a u(x)=-\displaystyle\lim_{y\rightarrow 0}y^av_y(x,y),
\end{equation*} and, for a given nonlinearity $ f\in C^{1,\gamma}(\R) $, we want to study the equation
\begin{equation}\label{eqtrace}
\mathcal{L}_au(x)=f(u) \qquad \text{in}\,\,\R^n.
\end{equation}
As a technical remark, we notice that,
in order to have the operator $ \mathcal{L}_a $ well defined for every smooth
function $ u:\R^n\to\R $, we need to choose the extension $ v $ in (\ref{lvsis}) in a unique
way. Indeed, for example, if~$v$ is a solution of~\eqref{lvsis} with~$a=0$,
then so is the function~$v(x,y)+e^{\pi x/2}\sin(\pi y/2)$.
To overcome this problem and uniquely determine~$v$
in~\eqref{lvsis}, we choose among all the possible solutions of (\ref{lvsis}) the one which is a minimizer of the energy
\begin{equation}\label{XCO}
\mathscr{E}_a(w):=\int_{\R^n\times(0,1)}y^a\abs{\nabla w(x,y)}^2\,dx\,dy,
\end{equation}
in the class of all the functions $w\in W^{1,2}(\R^n\times(0,1),y^a)$ such that $w(x,0)=u(x)$. Such a
minimizer $ v $ exists, it is unique, due to the convexity
of the energy functional in~\eqref{XCO},
and it solves the problem in~(\ref{lvsis}) --- see \cite{MV} for all the details.



With the setting
in~\eqref{lvsis}, the problem in~(\ref{eqtrace}) can be formulated in the following way:
\begin{align}
\label{mainsis}
\begin{cases}
\mathrm{div}(y^a\nabla v)=0 \qquad &\text{in}\,\,\R^n\times (0,1),\\
v_y(x,1)=0 \qquad &\text{on}\,\,\R^n\times \{y=1\},\\
-\displaystyle\lim_{y\rightarrow 0}y^av_y=f(v)\qquad &\text{on}\,\,\R^n\times\{y=0\},
\end{cases}
\end{align}
where $ f\in C^{1,\gamma}(\R)$ with $ \gamma>0 $.

Problem \eqref{mainsis} has a variational structure, since solutions of~(\ref{mainsis})
correspond to critical points of the energy functional
\begin{equation}\label{CRCR0}
\mathcal{E}(v):=\frac{1}{2}\int_{\R^n\times(0,1)}y^a |\nabla v (x,y)|^2 \,dx\,dy
+ \int_{\R^n\times \{y=0\}} F(v(x,0))\,dx,
\end{equation}
where the associated potential $F$ is such that $F'=-f$. 

Since problem (\ref{mainsis})
is set in a slab of fixed height, it is technically convenient to localize the energy
functional on cylinders. Namely, we define the cylinder
\begin{equation}\label{CRCR}C_R:=B_R\times(0,1),\end{equation} where $B_R\subset \R^n$
denotes the ball of radius $R$ centered at $0$. Then, by~\eqref{CRCR0},
the localized energy functional associated to problem~\eqref{mainsis} reads
\[
\mathcal{E}_R(v):=\frac{1}{2}\int_{C_R}y^a|\nabla v(x,y)|^2 \,dx\,dy
+ \int_{ B_R\times \{y=0\}} F(v(x,0))\,dx.
\]
In particular, the potential~$F$ is naturally defined up to an additive constant, hence,
focusing on bounded solutions,
we can also suppose that~$F\ge0$. For this kind of problems, the model case is the nonlinearity $ f(t):=t-t^3 $,
which arises in the study of phase transitions and it
is the derivative of the double-well potential
\[
F(t)=\frac14\left(1-t^2\right)^2.
\]

The usual notions of minimizer of the energy and of stable solution to problem (\ref{mainsis}) can be defined in a standard way.
We say that a bounded function $v\in C^1( \R^n\times (0,1))$  is a  \textit{minimizer} for~\eqref{mainsis} if
\[
\mathcal{E}_R(v)\leq \mathcal{E}_R(w)
\]
for every $R>0$ and for every bounded competitor $w$ such that $v\equiv w$ on $\partial B_R\times (0,1)$.

We say that a bounded solution $v$ of~\eqref{mainsis} is \textit{stable} if the second variation of the energy is non-negative, i.e.
\[
\int_{\R^n\times[0,1]}y^a\lvert \nabla \xi \rvert^2\,dx\,dy-\int_{\R^n\times \{y=0\}}f'(u)\xi^2\,dx \geq0
\]
for every function $\xi\in C_0^1(\R^n\times [0,1])$.

Clearly, if $v$ is a minimizer for~\eqref{mainsis} then, in particular, it is a stable solution. Another important subclass of stable solutions that we consider in this paper is given by the monotone solutions of (\ref{mainsis}). We say that a solution $ v $ of (\ref{mainsis}) is monotone if it is strictly monotone in one horizontal direction, say $ \partial_{x_n}v>0 $. For this kind of problems, it is possible to prove that monotone solutions are stable using a non-variational characterization of stability --- see for example Lemma 3.1 in \cite{CMV} for
all the details. See also~\cite{MR2779463} for a complete introduction to stable solutions in elliptic PDEs.

Problem (\ref{mainsis}) was initially studied by de la Llave and the third author in~\cite{DllV} with constant density, so with $ a=0 $. In particular, they proved a Liouville theorem that assures
the one-dimensional symmetry of monotone solutions on the trace,
provided that a suitable energy estimate for the functional associated to the problem
holds true. Since this energy estimate in dimension $ n=2 $ is a direct consequence of a
classical gradient bound, they obtain that monotone solutions of (\ref{mainsis}) with $ a=0 $ depend on only one horizontal variable if $ n=2 $. \medskip

We now describe some recent symmetry and rigidity results for problem~(\ref{mainsis})
in the light of a long-lasting line of investigation that was opened
by a celebrated conjecture by Ennio De Giorgi.

\subsection{Symmetry properties for the Allen-Cahn equation}
One of the main interests in proving the one-dimensional symmetry of monotone solutions comes from a
conjecture formulated by Ennio De Giorgi for the classical Allen-Cahn equation. Indeed, in 1979 De Giorgi posed the following
question. 

\begin{conjecture}\label{degiorgi}
	Let $u$ be a bounded and smooth solution of the Allen-Cahn equation
	\[
	-\Delta u=u-u^3 \qquad \text{in}\,\,\R^n,
	\]
	such that $ \partial_{x_n}u>0 $. Is it true that, if $ n\leq8 $, then $u$ is one-dimensional?
\end{conjecture}

A heuristic motivation of the conjecture can be formulated in light
of the work of Modica and Mortola \cite{MM}. Indeed, they proved that a proper rescaling
of the energy functional associated to the Allen-Cahn equation $ \Gamma $-converges to the perimeter functional, as the
rescaling parameter goes to zero. This means that a proper rescaling of the
minimizers of the Allen-Cahn equation converges to characteristic functions of sets
of minimal perimeter. The threshold dimension $ n=8 $ comes from the fact that
super-level sets of monotone functions are expected to be epigraphs
(though this is a tricky point, see e.g. formula~(5) in~\cite{FV}),
and minimal graphs are flat if $ n\leq8 $. For a complete discussion of minimal surfaces, see the illuminating monograph \cite{G}.

Summing up, the above heuristic argument would give that, at least in dimension~$n\le8$, if
we look at monotone solutions ``from very far'' (through a rescaling), their level sets are
close to hyperplanes. The question in Conjecture \ref{degiorgi} asks if, for this to hold,
the level sets of the function must be necessarily parallel hyperplanes.

The conjecture of De Giorgi
remained unanswered in every dimension $ n $ for almost twenty years.
It was proved to hold if $ n=2 $ by Ghoussoub and Gui \cite{GG} and by 
Berestycki, Caffarelli and Nirenberg~\cite{MR1655510},
and if $ n=3 $ by Ambrosio and Cabr\'e \cite{AC}.
Regarding dimensions $ 4\leq n\leq8 $, Savin
proved in~\cite{S} the conjecture by
assuming the following additional hypothesis about the limits in the monotone direction
\begin{equation}\label{ass_limits}
\lim_{x_n\to\pm\infty}u(x',x_n)=\pm1.
\end{equation}
Condition~\eqref{ass_limits} can be weakened by assuming two-dimensional
symmetry of the profiles at infinity, see~\cite{MR2728579}.

As a counterpart of the results giving positive answers to Conjecture~\ref{degiorgi}
(possibly under additional assumptions),
del Pino, Kowalczyk and Wei provided in \cite{dPKW} an example of a monotone solution to the Allen-Cahn equation in dimension $ n=9 $ which is not one-dimensional. In this way, they proved that dimension $ n=8 $ in Conjecture \ref{degiorgi} is the optimal one.

We refer to \cite{FV,CW} for more detailed surveys on topics related
to Conjecture~\ref{degiorgi}.
\medskip

\subsection{Symmetry properties for the fractional Allen-Cahn equation}
The fractional analogue of
Conjecture~\ref{degiorgi} can be formulated as follows:

\begin{conjecture}\label{frac_degiorgi}
	Let $s\in(0,1)$ and~$u$ be a bounded and smooth solution of the fractional Allen-Cahn equation
	\begin{equation}\label{frac_ac}
	(-\Delta)^s u=u-u^3 \qquad \text{in}\,\,\R^n,
	\end{equation}
	such that $ \partial_{x_n}u>0 $. Is it true that, if $ n$ is sufficiently small,
	then $u$ is one-dimensional?
\end{conjecture}

This question is also motivated by an analogue in the fractional setting of
the $ \Gamma $-convergence result by Modica and Mortola. Indeed, the third author
and Savin proved in \cite{SV} that a proper rescaling of the energy associated
to (\ref{frac_ac}) $\Gamma$-converges to the classical perimeter if~$ s\in\left[\frac12,1\right) $
and to the
fractional perimeter if~$ s\in\left(0,\frac12 \right)$. 

The fractional perimeter was introduced by Caffarelli, Roquejoffre and Savin in \cite{CRS},
and --- without going into the details --- can be thought as a nonlocal
version of the classical perimeter, counting the interactions between
points which lie in the two separated sides of the boundary of the set. As in the classical
case, one could relate, at least at a level of motivations, the validity
of Conjecture~\ref{frac_degiorgi} to the regularity and rigidity properties
of the minimizers of the limit energy functional, namely to the classical
minimal surfaces when~$s\in\left[ \frac12,1\right)$,
and to the nonlocal minimal surfaces when~$s\in\left(0,\frac12\right)$.
With respect to this, we recall that nonlocal minimal surfaces
are known to be smooth only in dimension~$2$ --- see~\cite{MR3090533} ---
and up to dimension~$7$ provided that~$s\in\left[\frac{1}{2}-\epsilon_0,\frac{1}{2} \right)$ and~$\epsilon_0$
is sufficiently small --- see~\cite{MR3107529}. Nonlocal minimal surfaces
that are entire graphs are known to be necessarily hyperplanes
only in dimension~$2$ and~$3$, and up to
up to dimension~$8$ provided that~$s\in\left[\frac{1}{2}-\epsilon_0,\frac{1}{2} \right)$ and~$\epsilon_0$
is sufficiently small --- see~\cite{MR3680376}.
Till now, no singular minimal surface is known --- 
see however~\cite{MR3798717} for the construction of a singular
cone in dimension~$7$ which is a stable critical point of the fractional
perimeter when~$s$ is sufficiently small.

Of course,
this lack of knowledge for the nonlocal minimal surfaces (when compared
to the classical minimal surfaces) provides a series of conceptual
difficulties when dealing with Conjecture~\ref{frac_degiorgi}, especially in the regime~$s\in\left(0,\frac12\right)$.

The problem posed by Conjecture \ref{frac_degiorgi} was solved in
dimension $n=2$ by Cabr\'e and Sol\`a-Morales in~\cite{CSM} for $s=\frac{1}{2}$,
and then by Cabr\'e, Sire and the third author in~\cite{YV, CY2}
for every~$s\in(0,1)$.

A positive answer in dimension $n=3$
was given by Cabr\'{e} and Cinti in~\cite{CC1}
and~\cite{CC2} in the cases~$s=\frac{1}{2}$
and~$s\in\left(\frac{1}{2},1\right)$, respectively. 
Regarding the strongly nonlocal regime, namely when~$s\in(0,\frac{1}{2})$,
recently the conjecture has been proved in dimension $n=3$
by Farina and the first and the third authors in~\cite{DFV} (using an improvement
of flatness result by~\cite{XFAH}) and by Cabr\'{e}, Cinti and Serra in~\cite{CCS}
(by a different approach which relies on some sharp energy
estimates and a blow-down convergence result for stable solutions).

Very recently, Figalli and Serra proved in~\cite{FS}
Conjecture \ref{frac_degiorgi} to be true for~$s=\frac12$ and~$ n=4 $
(also providing 
one-dimensional symmetry of stable solutions in dimension $ n=3 $).

Concerning higher dimensions, Savin proved in~\cite{S1,S2}
the conjecture for~$4 \leq n\leq 8$ and~$s\in\left[\frac{1}{2},1\right)$
under the additional
assumption~\eqref{ass_limits}.
Moreover, in~\cite{XFAH} it has been proved that 
Conjecture~\ref{frac_degiorgi} is true in dimensions~$4\le n\le 8$
if~$s\in\left[\frac12-\epsilon_0, \frac12\right)$, for some~$\epsilon_0$
sufficiently small, under the additional
assumption~\eqref{ass_limits}.

Besides these results,
Conjecture~\ref{frac_degiorgi} is also open in its generality, and the critical
dimension might depend on the fractional parameter~$s$. 

\subsection{Symmetry properties for the water wave problem}
Since, in our framework, we are dealing with a generalization
of fractional Laplace operators, which are attained in the limit according to~\eqref{ST},
a natural counterpart of Conjecture~\ref{frac_degiorgi} is the following one:

\begin{conjecture}\label{degiorgi-w}
	Let $a\in(-1,1)$ and~$u$ be a bounded and smooth solution of the fractional Allen-Cahn equation
	\[
	{\mathcal{L}}_a u=u-u^3 \qquad \text{in}\,\,\R^n,
	\]
	such that $ \partial_{x_n}u>0 $. Is it true that, if $ n$ is sufficiently small,
	then $u$ is one-dimensional?
\end{conjecture}

Conjecture~\ref{degiorgi-w} is related to, but structurally different from,
Conjecture~\ref{frac_degiorgi}. As a matter of fact, to point out the differences between problem (\ref{eqtrace})-(\ref{mainsis}) treated in these notes and its analogue for the fractional Laplacian, we consider the Fourier transform of the Dirichlet to Neumann operator $ \mathcal{L}_a $. It can be computed as
\begin{equation*}
\begin{aligned}
\widehat{ \mathcal{L}_a(u)}(\xi)=c_1(s)\frac{J_{1-s}(-i\abs{\xi})}{\cos\left(s\pi\right)
J_{1-s}(-i\abs{\xi})+\sin\left(s\pi\right) Y_{1-s}(-i\abs{\xi})}\abs{\xi}^{2s}\widehat{u}(\xi),
\end{aligned}
\end{equation*}
where $ J_m(\cdot) $ and $ Y_m(\cdot) $ are the Bessel functions of order $ m $, respectively of the first and second kind,
and~$ c_1(s) $ is a constant depending only on $ s\in(0,1) $.
As customary, the symbol~$\widehat{u}$ denotes the Fourier transform of~$u$.

Therefore, the operator $ \mathcal{L}_a $ can be seen as a Fourier operator with symbol
\[
S_s(\xi):=c_1(s)\frac{J_{1-s}(-i\abs{\xi})}{\cos\left(s\pi\right) J_{1-s}(-i\abs{\xi})+\sin\left(s\pi\right) Y_{1-s}(-i\abs{\xi})}\abs{\xi}^{2s}.
\]
The symbol $ S_s(\xi) $ was already known in \cite{BV,DllV} in the special case $ s=\frac12 $ as
\[
S_{1/2}(\xi)=\frac{e^{\abs{\xi}}-e^{-\abs{\xi}}}{e^{\abs{\xi}}+e^{-\abs{\xi}}}\abs{\xi}
\]
and it has been computed later by the second and third author in \cite{MV} for
every fractional parameter~$ s\in(0,1) $. By evaluating the limits of $ S_s(\xi) $ as $ \abs\xi $ goes to zero and infinity, we observe that
\begin{equation}\label{asympt}
\begin{aligned}
&S_s(\xi)\sim\abs{\xi}^2 \qquad &\text{as}\,\,\abs\xi\to0; \\
&S_s(\xi)\sim\abs{\xi}^{2s} \qquad &\text{as}\,\,\abs\xi\to+\infty.
\end{aligned}
\end{equation}

This fact is already evident in the simpler case $ s=\frac12 $, but it can
be shown also in the general case~$ s\in(0,1)$ --- see again \cite{MV} for all the details. To better undestand the implications of this behaviour, we should remind that $ \abs\xi^2 $ is the symbol of the classical Laplacian, and
that the fractional Laplacian can be also written in the Fourier setting as
\[
\widehat{ (-\Delta)^s u}(\xi)=\abs{\xi}^{2s}\widehat{u}(\xi),
\]
see for example \cite{H}. Looking at the asymptotics (\ref{asympt}), it becomes evident that the operator $ \mathcal{L}_a $ is not of purely fractional type, and, in fact,
it {\em
shows a nonlocal behaviour for high frequencies but it becomes similar to the Laplacian for small frequencies}.\medskip

In this setting, Conjecture~\ref{degiorgi-w} was first addressed
by de la Llave and the third author in \cite{DllV}, for the special case $ a=0 $. As mentioned above, their
main result is a Liouville theorem, that gives one-dimensional symmetry of
monotone solutions under an assumption about the growth of the Dirichlet energy
of the solution. In this way, they establish Conjecture~\ref{degiorgi-w} for~$a=0$ and~$n=2$
--- see in particular Theorem 1 in \cite{DllV}.

The results in~\cite{DllV} have been extended in \cite{CMV} from monotone to stable solutions,
also considering all the fractional parameters $ a\in(-1,1) $ and not only $ a=0 $.
In this setting, the result in \cite{CMV} reads as follows.

\begin{theorem}\label{DllV-T}
	Let $f\in C^{1,\gamma}(\R)$, with $\gamma>\max \{0,-a\}$, and let $v$ be a bounded and stable solution of~(\ref{mainsis}).
	
	Suppose that there exists $C>0$ such that
	\begin{equation}\label{EB}
	\int_{C_R}y^a |\nabla_x v(x,y)|^2\,dx\,dy\le
	CR^2
	\end{equation}
	for any $R\ge 2$, where the notation in~\eqref{CRCR} has been used.
	
	Then, there exist $v_0:\R\times(0,1)\rightarrow\R$ and $\omega
	\in {\rm S}^{n-1}$ such that
	\begin{equation*}
	v(x,y)=v_0 (\omega\cdot x,y)\qquad{\mbox{
			for any $(x,y)\in\R^{n}\times(0,1)$.}}
	\end{equation*}		
	In particular, the trace $u$ of $v$ on $\{y=0\}$ can be written as $u(x)=u_0(\omega\cdot x)$.
	
	Finally, either~$u_0'>0$ or $u_0'\equiv0$.
\end{theorem}

\begin{remark}\label{n=2}
{\rm	For this kind of elliptic problems, it is a standard
fact that bounded solutions have bounded gradients, see for example \cite{GT}.
For this reason, if we assume $ n=2 $, then hypothesis (\ref{EB}) is trivially
verified by any bounded stable solution. Therefore, we deduce that bounded stable solutions of \eqref{mainsis} are one-dimensional on the trace if $ n=2 $. In particular, this implies the
validity of Conjecture~\ref{degiorgi-w}
in dimension~$ n=2 $, for all~$a\in(-1,1)$ as a corollary of Theorem \ref{DllV-T}.}
\end{remark}

In~\cite{CMV}, Conjecture~\ref{degiorgi-w} is also addressed when $ n=3 $.
For this,
the strategy is based on energy estimates and the use of Theorem \ref{DllV-T}.
Namely, in \cite{CMV} the following result is proved:

\begin{theorem}[\textbf{Energy estimate for minimizers}]
	\label{globmin}
	Let $f\in C^{1,\gamma}(\R)$, with $\gamma>\max \{0,-a\}$, and let $v$ be a bounded minimizer for problem~\eqref{mainsis}. 
	
	Then, we have
	\begin{equation}
	\label{estimate.min}
	\mathcal{E}_R(v)=\frac{1}{2}\int_{C_R}y^a\lvert \nabla v \rvert^2 \,dx\,dy + \int_{ B_R\times \{y=0\}} F(v)\,dx 
	\leq C R^{n-1},
	\end{equation}
	for any $R\ge 2$, where the notation in~\eqref{CRCR} has been used.
	
\end{theorem}	
We point out that (\ref{estimate.min}) holds in general
for minimizers of the energy associated to problem~(\ref{mainsis}) in every dimension~$ n $,
but the application to symmetry problems usually becomes relevant only in dimension~$n\le3$.
Let us give now a brief look at the proof of Theorem \ref{globmin}, which is based on a direct comparison. Indeed, since we are assuming that $ v $ is a minimizer, for every admissible competitor $ w $ it holds that
\[
\mathcal{E}_R(v)\leq\mathcal{E}_R(w).
\] 
We say that a competitor $ w $ is admissible if $ w\equiv v $ on the lateral boundary $ \partial B_R\times(0,1)$. The key point for the proof is defining a competitor $ \widetilde{w} $ constantly equal to the minimum of the potential $ F $ in a cylinder of radius $ R-1 $, and then cutting it off in order to make it admissible. In such a way,
one is able to estimate the energy of $ \widetilde{w} $ in a cylinder of radius $ R $ and obtain (\ref{estimate.min}).
A strategy of this type has been used also in~\cite{AC} to solve the classical De Giorgi conjecture in
dimension~$ n=3 $.

Restricting to the case $n=3$, it is possible to prove the same estimate in Theorem~\ref{globmin}
for bounded solutions whose traces on $\{y=0\}$ are monotone in some direction, according to the following result:

\begin{theorem}[\textbf{Energy estimate for monotone solutions for $n=3$}]
	\label{thmono}
	Let $f\in C^{1,\gamma}(\R)$, with $\gamma>\max \{0,-a\}$, and let $v$ be a bounded solution of~(\ref{mainsis}) with $n=3$ such that its trace $u(x)=v(x, 0)$ is monotone in some direction. 
	
	Then, we have
	\begin{equation*}
	\mathcal{E}_R(v)=\frac{1}{2}\int_{C_R}y^a\lvert \nabla v \rvert^2 \,dx\,dy + \int_{ B_R\times \{y=0\}} F(v)\,dx 
	\leq C R^2,
	\end{equation*}
	for any $R\ge 2$, where the notation in~\eqref{CRCR} has been used.
\end{theorem}

We stress the fact that this energy estimate holds for monotone solutions of (\ref{mainsis}) only if we are in the case $ n=3 $. As mentioned above, this is due to the proof, and in particular to the fact that we know from Remark \ref{n=2} that stable solutions enjoy rigidity properties when we are in
one dimension less, i.e. when~$ n=2$.

Let us briefly sketch the proof of Theorem~\ref{thmono}. The interested reader can find all the details in Section~5 of \cite{CMV}.
First, it is necessary to define the two limit profiles of the monotone solution $ v $. Indeed, since $ v $ is monotone in one direction, say $ v_{x_3}>0 $, we can define
\[
\underline{v}(x_1.x_2,y):=\lim_{x_3\to-\infty}v(x,y);
\]
\[
\overline{v}(x_1,x_2,y):=\lim_{x_3\to+\infty}v(x,y).
\]
These functions are well defined for the monotonicity hypothesis, they are solutions of (\ref{mainsis}) in
one dimension less, so with $ n=2 $, and in particular they are stable. Here the
dimension plays a key role, since we can use Theorem \ref{DllV-T} and deduce
that $ \underline{v} $ and $ \overline{v} $ are
one-dimensional functions on the trace. From the existence of such solutions, it is possible to characterize the potential $F$. This is something that is fundamental in the proof.

On the other side, a monotone solution $ v $ is a minimizer of the energy in the class
\[
\mathcal{A}_v:=\left\{w\in H^1(C_R; y^a)\,\,\text{such that}\,\,w = v\,\,\text{on}\,\,\partial B_R\times(0,1)\,\,\text{and}\,\,\underline{v}\leq w\leq\overline{v}\,\,\text{in}\,\,C_R\right\}.
\]
For the detailed proof of this fact in the setting of water waves, see Lemma 5.6 in \cite{CMV}.

At this point, one can use the characterization of the potential $F$ provided by the previous steps in order to show that the competitor used in the proof of Theorem \ref{globmin} belongs to the class $ \mathcal{A}_v $, i.e. $ \widetilde{w}\in\mathcal{A}_v $.
Since the energy of $ \widetilde{w} $ in $ C_R $ can be bounded by~$ CR^2 $, this finishes the proof
of Theorem~\ref{thmono}.\medskip

The energy estimates in Theorems~\ref{globmin} and~\ref{thmono} give as a corollary the
following rigidity result for minimizers and for monotone solutions in dimension $ n=3 $.
Indeed, in this case hypothesis~(\ref{EB}) of Theorem~\ref{DllV-T} is fulfilled and the application is straightforward.

\begin{corollary}
	\label{1D}
	Let $f\in C^{1,\gamma}(\R)$, with $\gamma>\max \{0,-a\}$ and let $n=3$. Assume that one of the two following condition is satisfied:
	\begin{itemize}
		\item $v$ is a bounded minimizer for problem \eqref{mainsis};
		\item $v$ is a bounded solution of ~(\ref{mainsis}) such that its trace $u(x)=v(x,0)$ is monotone in some direction.
	\end{itemize}
	Then, there exist $v_0:\R\times(0,1)\to\R$ and $\omega\in S^2$ such that:
	\[
	v(x,y)=v_0(\omega\cdot x,y) \qquad {\mbox{for all }} (x,y)\in\R^3\times(0,1).
	\]
	
	In particular, the trace $u$ of $v$ on $\{y=0\}$ can be written as $u(x)=u_0(\omega\cdot x)$.
\end{corollary}

In particular, Corollary~\ref{1D} establishes the validity of
Conjecture~\ref{degiorgi-w} when $ n=3 $. The case of dimension~$n\ge4$ remains open.

\section*{Acknowledgement}

The first author has been supported by the DECRA Project DE180100957
``PDEs, free boundaries and applications''.
The first and third authors have been supported
by the Australian Research Council Discovery Project
DP170104880
``N.E.W. Nonlocal Equations at Work''.
The second author has been supported by MINECO
grant MTM2017-84214-C2-1-P and is part of the
Catalan research group 2017 SGR 1392. Part of this work was carried out on
the occasion of a very pleasant visit
of the second author to the University of Western Australia, which we
thank for the warm hospitality.

\end{document}